\documentclass[11pt,reqno]{amsart}
\usepackage{graphicx,amsmath,amsthm,verbatim,amssymb,lineno,titletoc}
\usepackage{mathrsfs}
\usepackage{adjustbox}
\usepackage{ytableau}
\usepackage{bbm}
\usepackage{hyperref}
\hypersetup{colorlinks}
\usepackage[numbers]{natbib}
\usepackage{array}
\usepackage{tikz}
\usetikzlibrary{shapes.multipart,patterns,arrows}

\usepackage[font=small]{caption}

\usepackage{amssymb}
\usepackage[T1]{fontenc} 
\usepackage{fourier} 
\usepackage[english]{babel} 
\usepackage{amsmath,amsfonts,amsthm,bbm} 
\usepackage{appendix}

\usepackage{amssymb}
\usepackage[all]{xy}
\usepackage{enumerate}
\usepackage{tikz}
\usepackage{mathrsfs}
\usepackage[english]{babel}
\usepackage{multirow}
\usepackage{float}
\usepackage{enumitem}
\usepackage{graphicx,accents}

\newcolumntype{M}[1]{>{\centering\arraybackslash}m{#1}}
\newcolumntype{N}{@{}m{0pt}@{}}

\pdfoutput=1

\newtheorem{innercustomthm}{Theorem}
\newenvironment{customthm}[1]
{\renewcommand\theinnercustomthm{#1}\innercustomthm}
{\endinnercustomthm}

\newtheorem{innercustomlemma}{Lemma}
\newenvironment{customlemma}[1]
{\renewcommand\theinnercustomlemma{#1}\innercustomlemma}
{\endinnercustomlemma}

\newtheorem{innercustomcorollary}{Corollary}
\newenvironment{customcorollary}[1]
{\renewcommand\theinnercustomcorollary{#1}\innercustomcorollary}
{\endinnercustomcorollary}

\newtheorem{innercustomprop}{Proposition}
\newenvironment{customprop}[1]
{\renewcommand\theinnercustomprop{#1}\innercustomprop}
{\endinnercustomprop}

\theoremstyle{definition}

\newtheorem{innercustomconjecture}{Conjecture}
\newenvironment{customconjecture}[1]
{\renewcommand\theinnercustomconjecture{#1}\innercustomconjecture}
{\endinnercustomconjecture}

\allowdisplaybreaks

\begin{document}
	
	\title[Chromatic number, induced cycles, and non-separating cycles]{Chromatic number, induced cycles, and non-separating cycles}

	\author{Hanbaek Lyu}
	\address{Hanbaek Lyu, Department of Mathematics, University of California Los Angeles, CA 90095}
	\email{\texttt{colourgraph@gmail.com}}

	\date{\today}
	
	\keywords{Chromatic number, induced cycles, Hadwiger number, induced non-separating cycles, Euler characteristic}

	\begin{abstract}
		We study two parameters obtained from the Euler characteristic by replacing the number of faces with that of induced and induced non-separating cycles. By establishing monotonicity of such parameters under certain homomorphism and edge contraction, we obtain new upper bounds on the chromatic number in terms of the number of induced cycles and the Hadwiger number in terms of the number of induced non-separating cycles. As an application, we show that a 3-connected graph with average degree $k\ge 2$ have at least $(k-1)|V|+Ck^{3}\log^{3/2}k$ induced non-separating cycles for some explicit constant $C>0$. This improves the previous best lower bound $(k-1)|V|+1$, which follows from Tutte's cycle space theorem. We also give a short proof of this theorem of Tutte. 
	\end{abstract}
	
	${}$
	\vspace{-0.5cm}
	${}$
	\maketitle

	\section{Introduction}
	\label{Introduction}

	Understanding structure of graphs is often a formidable task. But sometimes, simply counting the number of basic objects such as vertices, edges, and cycles, gives a good understanding on some structural properties of the graphs. This is best illustrated by a classic theorem of Euler, which says that the alternating sum of the number of faces, edges and vertices in a graph $G$ determines its genus. Namely, if a graph $G$ is obtained by a 2-cell division of an orientable surface of genus $g$, then    
	\begin{equation}\label{eq:euler_characteristic}
	|\mathfrak{F}(G)|-|E(G)|+|V(G)|=2-2g
	\end{equation}  
	where $\mathfrak{F}(G)$ is the set of all faces in $G$. The quantity in the left hand side is called the \textit{Euler characteristic}.
	
	The notion of faces depends not only on the graph itself but also on the underlying surface on which the graph is drawn. However, it is well known that face boundaries of a 3-connected planar graphs are precisely its induced non-separating  cycles. Hence, it is natural to replace the set $\mathfrak{F}(G)$ of all faces in the Euler characteristic by some other class of cycles $\mathcal{C}(G)$. Parameters obtained in this way could be used to study different structural properties of graphs other than their genus. Indeed, Vince and Littel \cite{vince1989discrete} and Yu \cite{yu1992non} used the parameter obtained when $\mathcal{C}(G)$ is a cycle double cover (a collection of cycles such that each edge is contained in exactly two cycles in $\mathcal{C}(G)$) in oder to study discrete Jordan curves.  
	
	In this paper, we study two parameters obtained when $\mathcal{C}(G)$ is the set $C(G)$ of all induced cycles and the set $F(G)$ of all induced non-separating cycles in $G$. Let $\chi(G)$ denote the chromatic number of $G$. Our first main result is the following.
	
	\begin{customthm}{1}\label{thm:chromatic}
		Let $G=(V,E)$ be a connected graph. Then we have 
		\begin{equation}\label{eq:thm1}
		\binom{\chi(G)}{3} - \binom{\chi(G)}{2}+\binom{\chi(G)}{1} \le |C(G)|-|E(G)|+|V(G)|.
		\end{equation}
	\end{customthm}
	Note that Theorem \ref{thm:chromatic} gives an upper bound on $\chi(G)$ in the order of $|C(G)|^{1/3}$. This is incomparable with the easy bound $\chi(G)\le 3k+2$, where $k$ is the maximum number of vertex-disjoint induced cycles in $G$.     
	
	Let  $h(G)$ denote the Hadwiger number of $G$, the size of maximum complete graph minor in $G$. Our second main result is the following, which gives an upper bound on $h(G)$ in the order of $|F(G)|^{1/3}$ when $G$ is 3-connected. 
	
	\begin{customthm}{2}\label{thm:hadwiger}
		Let $G=(V,E)$ be a 3-connected graph. Then we have 
		\begin{equation}\label{eq:thm2}
		\binom{h(G)}{3} - \binom{h(G)}{2}+\binom{h(G)}{1} \le |F(G)|-|E(G)|+|V(G)|.
		\end{equation}
	\end{customthm}
	We prove Theorems \ref{thm:chromatic} and \ref{thm:hadwiger} in Sections \ref{section:pf_chromatic} and \ref{section:pf_minor}, respectively.

	As an application, Theorem \ref{thm:hadwiger} gives a new lower bound on the number of induced non-separating cycles in 3-connected graphs. Note that the previous best lower bound follows from a theorem of Tutte, which states that the cycle space of a 3-connected graph is generated by its induced non-separating cycles \cite{tutte1963draw}. Recall that the rank of the cycle space of a connected graph $G$ equals that of its fundamental group, which is $|E(G)|-|V(G)|+1$. Hence this gives 
	\begin{equation}\label{eq:F_lowerbd_tutte}
	|F(G)| \ge |E(G)| - |V(G)|+1.
	\end{equation} 
	In fact, this inequality follows directly from Theorem \ref{thm:hadwiger}, since the left hand side of (\ref{eq:thm2}) is at least 1.

	Graphs with large average degree are known to have large Hadwiger number due to Kostochka \cite{kostochka1984lower}. Namely, for a graph $G=(V,E)$ with average degree $|E|/|V|\ge k$ for a fixed integer $\kappa\ge 2$, we have $h(G)\ge k/270 \sqrt{\log k}$. Noting that the left hand side of (\ref{eq:thm2}) is non-decreasing in $h(G)$, Theorem \ref{thm:hadwiger} yields the following new lower bound on $|F(G)|$:

	\begin{customcorollary}{3}\label{lowerbd}
		Let $G=(V,E)$ be a 3-connected graph with $|E|/|V| \ge k$ for some integer $k \ge 2$. Then we have 
		\begin{equation}\label{avgdegineq}		
		|F(G)|\ge (k-1)|V|+\frac{1}{118098000} \left(k^{3}\log^{3/2}k - 1620k^{2}\log k + 801900k\log^{1/2} k\right) .
		\end{equation}
	\end{customcorollary}
	
	We remark that this improves the previous best lower bound (\ref{eq:F_lowerbd_tutte}), which reads $|F(G)|\ge (k-1)|V|+1$.

	\subsection{Definitions}
	\label{subsection:def}
	
	In this paper, every graph is finite and simple. Let $S,H\subseteq G$ be two subgraphs in $G$. Define $S+H$ to be the subgraph of $G$ with $V(S+H)=V(S)\cup V(H)$ and edge set $E(S+H)=E(S)\cup E(H) \cup \{ \text{edges between $S$ and $H$} \}$. We let $S-H:=S-V(H)$ denote the subgraph of $G$ obtained from $S$ by deleting all vertices of $H$ and edges incident to them.  If $H$ is the singleton $\{v\}$ for some vertex $v$ in $G$, then we denote $S+v$ and $S-v$ for $S+H$ and $S-H$, respectively. 

	If $H\subseteq G$ and $v\in V(H)$, then $N_{H}(v)$ denotes the set of all neighbors of $v$ in $H$ and $\deg_{H}(v):=|N_{H}(v)|$. A graph $G$ is a \textit{complete graph} if $uv\in E(G)$ for all distinct vertices $u,v$ in $G$. A complete graph with $n$ vertices is denoted by $K_{n}$. A graph $C$ is called a \textit{cycle} if it is connected and 2-regular, meaning $\deg_{C}(v)=2$ for all $v\in V(C)$. We denote by $C_{n}$ the cycle of $n$ vertices. A subgraph $H\subseteq G$ is \textit{induced} if $S\subseteq H$ for any subgraph $S \subseteq G$ such that $V(S)=V(H)$, and \textit{non-separating} if $G-H$ is connected. We let $C(G)$ and $F(G)$ denote the set of all induced and induced non-separting cycles in $G$, respectively. We denote the parameters in the right hand sides of (\ref{eq:thm1}) and (\ref{eq:thm2}) by $\Lambda_{C}(G)$ and $\Lambda_{F}(G)$, respectively. A graph $G$ is \textit{3-connected} if $|V(G)|\ge 4$ and $G-H$ is connected for all subgraph $H$ such that $|V(H)|\le 2$.

	\vspace{0.3cm}

\section{Proof of Theorem \ref{thm:chromatic}}
\label{section:pf_chromatic}

In this section, we prove Theorem \ref{thm:chromatic}. Given two graphs $G$ and $H$, a map $\nu:V(G)\rightarrow V(H)$ is called a \textit{homomorphism} if it preserves the adjacency relation, i.e., $uv\in E(G)$ implies $\nu(u)\nu(v)\in E(H)$. We denote $\nu:G\rightarrow H$ for a homomorphism $\nu:V(G)\rightarrow V(H)$. For example, the map from a graph to the graph obtained by identifying two nonadjacent vertices and deleting resulting loops and parallel edges is a homomorphism, sometimes called the \textit{nonedge contraction}. We say $\nu:G\rightarrow H$ is an \textit{isomorphism} if there exists a homomorphism $\xi:H\rightarrow G$ such that $\xi\circ\nu:G\rightarrow G$ is the identity map on $V(G)$. A \textit{n}-coloring of $G$ is a homomorphism $G\rightarrow K_{n}$. The \textit{chromatic number} of $G$, denoted by $\chi(G)$, is the least integer $n$ such that $G$ admits a $n$-coloring. 

\begin{figure*}[h]
	\centering
	\includegraphics[width=0.5 \linewidth]{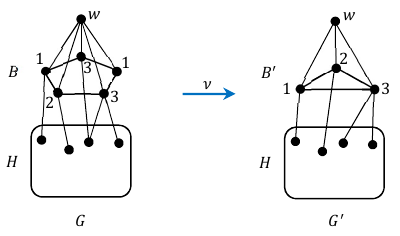}
	\vspace{-0.2cm}
	\caption{ An example of local homomorphism $\nu:G\rightarrow G'\in [G/w]$, which maps $B=G[N(w)]\simeq C_{5}$ to $B'\simeq K_{3}$ according to the 3-coloring of $C_{5}$ shown above, while mapping $G-B$ isomorphically.
	}
	\label{fig:compression_ex}
\end{figure*}

We introduce an elementary homomorphism, which will be the basis for our proof of Theorem \ref{thm:chromatic}. Given a graph $G$ and a vertex $w\in V(G)$, let $B=G[N(w)]$ be the induced subgraph of $G$ on the vertex set $N(w)$. A \textit{local homomorphism} of $G$ at $w$ is a homomorphism $\nu:G\rightarrow G'$ such that $\nu$ restricted on $G-B$ is an isomorphism. Define a set of graphs $[G/w]$ by
\begin{equation}
[G/w] = \left\{ G' \,\,\bigg| \,\,  \begin{matrix} \text{$\nu:G\rightarrow G'$ is a local homomorphism of $G$ at $w$} \\ \text{and $\nu(B)$ is isomosphic to $K_{\chi(B)}$ }   \end{matrix} \right\}.
\end{equation}
In words, $[G/w]$ is the set of all graphs obtained by successively identifying pairs of vertices in $N(w)$ such that $B$ is mapped to the smallest complete graph $K_{\chi(B)}$ (see Figure \ref{fig:compression_ex}).

Recall that $\Lambda_{C}(G)$ denotes the parameter in the right hand side of (\ref{eq:thm1}).  The key observation in proving Theorem \ref{thm:chromatic} is the following monotonicity of the parameter $\Lambda_{C}$ under a local homomorphism $G\rightarrow G'\in [G/w]$, as stated in the following lemma.

\begin{customlemma}{2.1}\label{lemma:compression}
	Let $G=(V,E)$ be a connected graph and fix $w\in V$. Then for any $G'\in [G/w]$, we have 
	\begin{equation}\label{eq:lemma_compression}
	\Lambda_{C}(G')\le \Lambda_{C}(G).
	\end{equation}	 
\end{customlemma}

\hspace{-0.43cm}\textbf{Proof of Theorem \ref{thm:chromatic}.} If $G$ is not a complete graph, then we can choose a vertex $w\in V(G)$ such that for any $G'\in [G/w]$, $|V(G')|<|V(G)|$. Hence we may choose a sequence of graphs $G=G_{0}$, $G_{1},\cdots, G_{k}$ such that $G_{k}$ is a complete graph and for all $0\le i < k$ we have $G_{i+1}\in [G_{i}/w_{i}]$ for some $w_{i}\in V(G_{i})$. Then the composition of homomorphisms $f:G_{0}\rightarrow G_{1}\rightarrow \cdots \rightarrow G_{k}$ is a $|V(G_{k})|$-coloring on $G$ so we have $\chi(G)\le |V(G_{k})|$. By Lemma \ref{lemma:compression}, we have $\Lambda_{C}(G_{n})\le \Lambda_{C}(G)$. Note that 
\begin{equation}
\Lambda_{C}(K_{n}) = \binom{n}{3} - \binom{n}{2}+\binom{n}{1} 
\end{equation}
is a non-decreasing function in $n$. Hence $\Lambda_{C}(K_{\chi(G)})\le \Lambda_{C}(G_{k}) \le \Lambda_{C}(G)$, as desired. $\hfill\blacksquare$
\vspace{0.2cm}

In the rest of this section, we prove Lemma \ref{lemma:compression}. We begin by describing how induced cycles behave under local homomorphisms.

\begin{customlemma}{2.2}\label{lemma:cycleinjection}
	Let $G=(V,E)$ be a connected graph with a vertex $w$. Fix a local homomorphism $\nu:G\rightarrow G'\in [G/w]$. Denote $B=G[N(w)]$, $H=G-B-w$, and $B'=\nu(B)$. Then we have the following.
	\begin{description}
		\item[(i)] There exists an injection $\phi : C(G')\setminus C(B'+w) \rightarrow C(G)\setminus C(B+w)$ such that for all $C\in C(G')\setminus C(B'+w)$, $V(\nu(\phi(C))) \setminus \{w\}=V(C)$.
		\vspace{0.1cm}
		\item[(ii)] $|C(G)\setminus C(B+w)|- |C(G')\setminus C(B'+w)| \ge |E(G)\setminus E(B+w)|- |E(G')\setminus E(B'+w)|$.
	\end{description}
\end{customlemma}

\begin{proof}
	We first show (i). Recall that $\nu$ is an isomorphism between $G-B$ and $G'-B'$, so we may identify the two graphs by $\nu$. Fix an induced cycle $C$ in $G'$ that is not contained in $B'+w$. We shall correspond an induced cycle $\phi(C)$ in $G$ such that the property stated at the end of (i) holds. Observe that $C$ does not use the vertex $w$ and contains at most two vertices in $B'$, since otherwise it would be a triangle contained $B'+w$. If $C$ uses no vertex in $B'+w$, then it is also an induced cycle in $G$ so we define $\phi(C)=C$.

	Suppose $C$ uses one vertex, say $z$, in $B'$. Let $x,y$ be the two neighbors of $z$ in $C$. Then $P:=C-z\subseteq H$ is an induced path from $x$ to $y$ both in $G/w$ and $G$. If there exists a vertex $v$ in $B$ such that $\nu(v)=z$ and $vx,vy\in E(G)$, then define $\phi(C)=P+v$ (see Figure \ref{fig:induced_cycle1} (a)). If not, there exists two vertices $u,v$ in $B$ such that $\nu(u)=\nu(v)=z$ and $uv,vy\in E(G)$. Note that $uv\notin E(G)$ since the homomorphism $\nu$ identifies $u$ and $v$ into the single vertex $z$. So $P+u+v+w$ is an induced cycle in $G$, and we define this to be $\phi(C)$ (see Figure \ref{fig:induced_cycle1} (b)).
	
	\begin{figure*}[h]
		\centering
		\includegraphics[width=0.7 \linewidth]{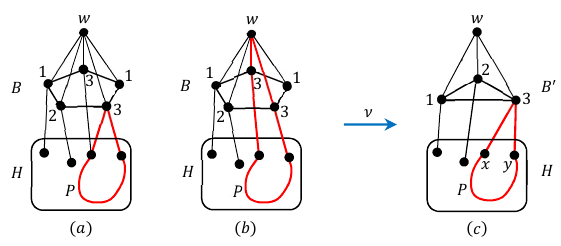}
		\vspace{-0.2cm}
		\caption{ Definition of $\psi(G)$ when $C\in C(G')\setminus C(B'+w)$ uses one vertex in $B'$. The red cycle in (c) is $C$, and $\phi(C)$ is the red cycles in (a) and (b) in each cases.
		}
		\label{fig:induced_cycle1}
	\end{figure*}
	
	Lastly, suppose $C$ uses two vertices, say $z_{1}$ and $z_{2}$, in $B'$. Let $x,y$ be the neighbors (not necessarily distinct) of $z_{1},z_{2}$ in $C$ repectively. Then $P:=C-z_{1}-z_{2}\subseteq H$ is an induced path from $x$ to $y$ both in $G'$ and $G$. Since $\nu$ is a homomorphism, there exists vertices $v_{1},v_{2}$ in $B$ such that $\nu(v_{1})=z_{1}$, $\nu(v_{2})=z_{2}$, and $xv_{1},yv_{2}\in E(G)$. Now if $v_{1}v_{2}\in E(G)$ we define $\phi(C)=P+v_{1}+v_{2}$ (see Figure \ref{fig:induced_cycle2} (a)) and otherwise define $\phi(C)=P+v_{1}+v_{2}+w$ (see Figure \ref{fig:induced_cycle2} (b)). This defines the map $\phi:C(G')\setminus C(B'+w)\rightarrow C(G)\setminus C(B+w)$ asserted in (i). From the construction it is clear that the last part of (i) holds. The injectivity of $\phi$ follows from this. This shows (i).
	
	\begin{figure*}[h]
		\centering
		\includegraphics[width=0.7 \linewidth]{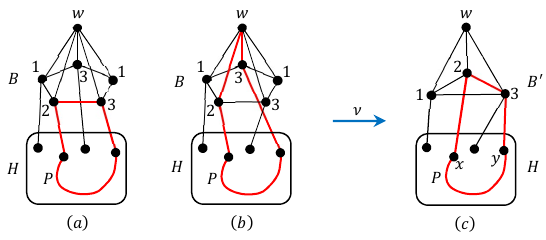}
		\vspace{-0.2cm}
		\caption{ Definition of $\psi(G)$ when $C\in C(G')\setminus C(B'+W)$ uses two vertices in $B$. The red cycle in (c) is $C$, and $\phi(C)$ is the red cycles in (a) and (b) in each cases.
		}
		\label{fig:induced_cycle2}
	\end{figure*}

	Next we show (ii). It suffices to show that the number of induced cycles in $C(G)\setminus C(B+w)$ that are not in the image of $\phi$ is at least the right hand side of the inequality in (ii). Let $R$ the set of ``wedges'' in $G-w$ that becomes a $K_{2}$ by $\nu$, i.e., 
	\begin{equation}
	R = \{ \text{$W$ a path in $G$}\,|\, \text{ $V(W)=\{x,y,z\}$, $\{x,z\}\subseteq V(B)$, $y\in V(H)$, and $\nu(x)=\nu(z)$} \}
	\end{equation} 
	Note that for each $W\in R$, $W+w$ is an induced $C_{4}$ in $G$ whose image under $\nu$ is a length 2 path in $G/w$. By the property of $\phi$, these rectangles are not in the image of $\phi$. On the other hand, each such wedge $W\in R$ corresponds to a pair of edges in $E(G)\setminus E(B+w)$ that are identified by $\nu$. Thus $|R|$ equals to the right hand side of (ii). This shows the assertion. 
\end{proof}

The following proposition will be useful in the proof of Theorem \ref{thm:chromatic}.

\begin{customprop}{2.3}\label{prop:chi_induction_step}
	Suppose Lemma \ref{lemma:compression} holds for every pair $(G,w)$ of connected graph $G$ and a vertex $w\in V(G)$ whenever $|V(G)|\le n$ for some integer $n\ge 1$. Then for any graph $H$ not necessarily connected and $|V(H)|<n$, we have  
	\begin{equation}
	|C(H)| \ge |C(K_{\chi(H)})|.
	\end{equation}
\end{customprop}

\begin{proof}
	We first observe that for any connected graph $S$ such that $|V(S)|<n$,
	we have
	\begin{equation}\label{eq:obs_pf_chromatic}
	|C(S)|\ge |C(S')|
	\end{equation}
	for any $u\in V(S)$ and for all $S'\in [S/u]$. Indeed, the assumption implies
	\begin{equation}
	|C(S)|-|C(S')| \ge (|E(S)|-|E(S')|) - (|V(S)|-|V(S')|) \ge 0,
	\end{equation}
	where the second inequality follows since the homomorphism $S\rightarrow S'\in [S/u]$ identifies more edges than vertices; for, if two non-adjacent vertices $x$ and $y$ in $N_{S}(u)$ are identified, then the two edges $xu$ and $yu$ are also identified. 
	
	Now, let $H_{1},\cdots,H_{k}$ be the components of $H$. For each $H_{i}$, choose a sequence of graphs $H_{i}=H_{i}(0),H_{i}(1),\cdots, H_{i}(r_{i})=H_{i}'$ such that $H_{i}'$ is a complete graph and we have ${H_{i}(j+1)}\in [H_{i}(j)/w_{i}(j)]$ for some $w_{i}(j)\in V(H_{i}(j))$ for all $0\le j < r_{i}$. Since $|V(H_{i})|\le |V(H)|<n$,  one can apply (\ref{eq:obs_pf_chromatic}) so that 
	\begin{equation}
	|C(H_{i})|\ge |C(H_{i}')| \qquad 1\le i \le k.
	\end{equation}
	
	Next, embed all complete graphs $H_{i}'$ into a maximal one, say $H_{*}'$, by some homomorphism. Note that 
	\begin{equation}\label{eq:prop_inductionstep}
	|C(H)| = \sum_{i=1}^{k}|C(H_{i})| \ge \max_{1\le i \le k}(|C(H_{i}')|) = |C(H_{*}')|.
	\end{equation}
	On the other hand, the composition of all homomorphisms we have used gives a coloring $H\rightarrow H_{*}'$. By the minimality of $\chi(H)$, we have $ |V(H_{*}')|\ge \chi(H)$. Since $H_{*}'$ is a complete graph, this yields $ |C(H_{*}')|\ge |C(K_{\chi(H)})|$. Then the assertion follows from (\ref{eq:prop_inductionstep}).  
\end{proof}

Now we are ready to give a proof of Lemma \ref{lemma:compression}.

\vspace{0.2cm}
\hspace{-0.425cm}\textbf{Proof of Lemma \ref{lemma:compression}.}  We use induction on $|V|$. We may assume $|V|>1$ since otherwise the assertion is trivial. For the induction step, let $n\ge 2$ and suppose for any pair $(G,w)$ of connected graph $G$ and a vertex $w\in V(G)$ that the assertion holds whenever $|V(G)|<n$. We use the same notation as in the statement of Lemma \ref{lemma:cycleinjection}. Denote $\Delta C:=|C(G)|-|C(G')|$,  $\Delta E:=|E(G)|-|E(G')|$, and $\Delta V:=|V(G)|-|V(G')|$.

Since $B'\simeq K_{\chi(B)}$ and $w\notin B$, by the induction hypothesis and Proposition \ref{prop:chi_induction_step}, we have 
\begin{equation}
|C(B)|\ge |C(B')|.
\end{equation}
It is easy to see that $|C(B+w)|=|C(B)|+|E(B)|$ and similarly for $B'+w$. Hence we have 
\begin{align}
|C(B+w)|-|C(B'+w)| &= |C(B)|-|C(B')| + |E(B)|-|E(B')| \\
&\ge |E(B)|-|E(B')|.
\end{align}
Then note that   
\begin{eqnarray}
|E(B+w)|-|E(B'+w)|&=& |E(B)| -|E(B')|+|V(B)|-|V(B')|\\
&=& |E(B)| -|E(B')|+\Delta V.
\end{eqnarray}
Thus by Lemma \ref{lemma:cycleinjection} (ii),  we have 
\begin{eqnarray}
\Delta C &=&  (|C(B+{w})|-|C(B'+w)|)+\left( |C(G)\setminus C(B+w)|- |C(G')\setminus C(B'+w)| \right)\\
&\ge & |E(B)|-|E(B')|+|E(G)\setminus E(B+w)|- |E(G')\setminus E(B'+w)|\\
&=& \Delta E - \Delta V,
\end{eqnarray}
as desired. This shows the assertion. $\hfill\blacksquare$

\section{Proof of Theorem \ref{thm:hadwiger}}
\label{section:pf_minor}

In this section, we prove Theorem \ref{thm:hadwiger}. Given a graph $G$ and an edge $e=uv\in E(G)$ let $G/e$ be the graph obtained by identifying vertices $u$ and $v$ into a single vertex $v_{e}$, where we automatically suppress any multiple edges and loops so that the resulting graphs are always simple. This operation $G\rightarrow G/e$ is called  \textit{edge contraction}. We say a graph $H$ is a \textit{minor} of $G$ if it can be obtained by a sequence of edge contractions from a subgraph $H'\subseteq G$. We write $H\le_{m} G$ if $H$ is a minor of $G$. The \textit{Hadwiger number} of $G$, denoted $h(G)$, is the maximum number $n$ such that $K_{n}\le_{m} G$. 

Let $\nu_{e}:V(G)\rightarrow V(G/e)$ be the vertex map induced by the edge contraction $G\rightarrow G/e$. Note that gives a homomorphism $G-e\rightarrow G/uv$, which identifies the non-adjacent vertices $u,v$ in $G-e$ into $v_{e}$ and suppresses any resulting loops and parallel edges. For any subgraph $H\subseteq G$, we denote by $H/e$ the subgraph $\nu_{e}(H-e)\subseteq G/e$. Note that this agrees with the usual edge contraction when $e\in E(H)$. If $u\in V(H)$ and $v\notin V(H)$, then $H/e$ equals the graph obtained from $H$ by renaming $u$ as $v_{e}$. If both $u$ and $v$ are not in $H$, then $H/e= H$.

Recall that $\Lambda_{F}(G)$ denotes the parameter in the right hand side of (\ref{eq:thm2}). The following is a key observation in proving Theorem \ref{thm:hadwiger}.

\begin{customlemma}{3.1}\label{lemma:minor}
	Let $G=(V,E)$ be a 3-connected graph and let $e\in E$ be such that $G/e$ is 3-connected. We have 
	\begin{equation}
	\Lambda_{F}(G/e)\le \Lambda_{F}(G).
	\end{equation}	 
\end{customlemma}

In order to prove Theorem \ref{thm:hadwiger}, we need the existence of an edge $e\in E(G)$ such that $G/e$ remains 3-connected and $h(G/e)=h(G)$, provided $G$ is not a complete graph. This is provided by Seymour's celebrated ``Splitter theorem'', which is stated below in a simple form.

\begin{customthm}{3.2}[Seymour \cite{seymour1986decomposition}]\label{seymour}
	Let $H\le_{m} G$ where both $H$ and $G$ are 3-connected and $|V(G)|>4$. If $|V(H)|<|V(G)|$, then there exists an edge $e\in E(G)$ such that $G/e$ is 3-connected and contains $H$ as a minor.
\end{customthm}

\vspace{0.2cm}
\hspace{-0.42cm}\textbf{Proof of Theorem \ref{thm:hadwiger}.} Let $G=(V,E)$ be a 3-connected graph. We use induction on $|V|$. Let $K$ be the largest complete graph that is a minor of $G$. So $K=K_{h(G)}$. Noting that $\Lambda_{F}(K_{n}) = \binom{n}{3} - \binom{n}{2}+\binom{n}{1}$, it suffices to show 
\begin{equation}
\Lambda_{F}(K)\le \Lambda_{F}(G).
\end{equation}
We may assume $K\ne G$ since otherwise there is nothing to prove.  If $|V|=4$, then 3-connectedness yields $G=K_{4}$ so the assertion holds. For the induction step, we may assume the assertion holds for 3-connected graphs with less than $|V|$ vertices. Since $K\subseteq G$ is a complete graph, $|V(K)|=|V|$ implies that $G$ is a complete graph so the assertion holds. We may assume $|V(K)|<|V(G)|$. Then by Theorem \ref{seymour}, there exists an edge $e$ in $G$ such that $G/e$ is 3-connected and still contains $K$ as a minor. Thus $K\le_{m} G/e \le_{m} G$, which yields $h(G/e)=h(G)$. Then by the induction hypothesis and Lemma \ref{lemma:minor}, we obtain 
\begin{equation*}
\Lambda_{F}(K)=\Lambda_{F}(K_{h(G/e)}) \le \Lambda_{F}(G/e) \le \Lambda_{F}(G),
\end{equation*} 
as desired. $\hfill\blacksquare$
\vspace{0.2cm}

Next, we prove Lemma \ref{lemma:minor}. To this end, we first need to understand how induced non-separating cycles in a 3-connected graph behave when contracting an edge $e$ such that the resulting graph is still 3-connected. 

\begin{figure}[H]
	\centering
	\includegraphics[width=0.9 \linewidth]{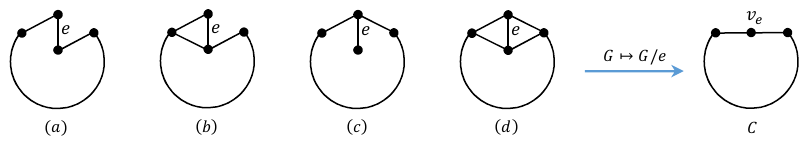}
	\caption{Possible types of inverse image $C*e$ of an induced cycle $C$ in $G/e$ using the contracted vertex $v_{e}$ under the edge contraction $G\rightarrow G/e$. }
	\label{fig:contractionmap}
\end{figure}

Let $G$ be a graph and fix an edge $e=uv\in E(G)$. For each $C\in C(G/e)$ such that $v_{e}\in V(C)$, note that $C-v_{e}$ is an induced path both in $G$ and $G/e$. We denote $C*e = (C-v_{e})+u+v\subseteq G$. Note that $C*e$ contains the edge $e$ and it is the largest subgraph (w.r.t. inclusion) in $G$ such that $(C*e)/e=C$. Moreover, observe that $C*e$ is one of the four types in Figure \ref{fig:contractionmap} up to isomorphism, depending on the degrees of vertices $u$ and $v$ in $C*e$. 

\begin{customprop}{3.3}\label{prop:injection1}
	Let $G=(V,E)$ be a 3-connected graph and $e=uv\in E$ be an edge such that $G/e$ is 3-connected. Suppose $C\in F(G/e)$ and $v_{e}\in V(C)$. Then the following holds. 
	\begin{description}
		\item[(i)] If $\deg_{C*e}(u)=\deg_{C*e}(u)=2$, then $C*e\in F(G)$. 
		\vspace{0.1cm}
		\item[(ii)] If $\deg_{C*e}(v)<3$, then $(C-v_{e})+u\in F(G)$. 
		\vspace{0.1cm}
		\item[(iii)] If $\deg_{C*e}(u)=\deg_{C*e}(u)=3$, then either $(C-v_{e})+u$ or $(C-v_{e})+v$ is in $F(G)$.
	\end{description}
\end{customprop}

\begin{proof}

	For (i), suppose $\deg_{C*e}(u)=\deg_{C*e}(u)=2$ (see Figure \ref{fig:contractionmap} (a)). Note that  ${C*e} = (C-v_{e})+u+v\in C(G)$ and $G-(C*e)=(G/e)-C$ is connected since $C\in F(G/e)$. Hence $C*e\in F(G)$. For (ii), suppose $\deg_{C*e}(v)<3$ (see Figure \ref{fig:contractionmap} (b) and (c)). Then $C':=(C-v_{e})+u\in C(G)$. Since $G-(C*e)=G/e - C$ is connected and $C*e = C'+v$, $G-C'$ is connected if $v$ is adjacent to some vertex in $G-(C*e)$. If not, since $C*e$ is an induced subgraph in $G$, it follows that $\deg_{G}(v)<3$. But this contradicts 3-connectivity of $G$. Hence $C'\in F(G)$.

	For (iii), note that both $C_{1}:=(C-v_{e})+u$ and $C_{2}=(C-v_{e})+v$ are induced cycles in $G$. Suppose, for a contradiction, that both $C_{1}$ and $C_{2}$ are separating. Let $B=G-(C*e)\subseteq G$, which is nonempty since $G\ne K_{4}$. If $u$ is adjacent to $B$ in $G$, then $C_{2}$ is non-separating. So we may assume $u$ is not adjacent to $B$. Similarly, we may assume $v$ is not adjacent to $B$. Let $L=u+v\subseteq G$. Then $G-x-y$ is not connected, where the nonempty subgraphs $L$ and $B$ in $G$ belong to different components. But this contradicts the 3-connectivityf of $G$. This shows (iii). 
\end{proof}

Based on the previous proposition, it is now straightforward to construct an injection $F(G/e)\rightarrow F(G)$ when both $G$ and $G/e$ are 3-connected. 

\begin{customprop}{3.4}\label{prop:injection2}
	Let $G=(V,E)$ and $e=uv\in E$ be such that both $G$ and $G/e$ are 3-connected. Then there exists an injection $\psi:F(G/e)\rightarrow F(G)$ such that $\psi(C)/e=C$ for all $C\in F(G/e)$. 
\end{customprop}

\begin{proof}
	Let $C\in F(G/e)$. If $C$ does not use $v_{e}$, then $C\subseteq G$ and  $(G-C)/e = G/e - C$ so $C\in F(G)$. Hence we define $\psi(C)=C$. Otherwise, we let $\psi(C)$ be one of $(C-v_{e})+u$, $(C-v_{e})+v$, or $(C-v_{e})+uv$, according to Proposition \ref{prop:injection1} so that $\psi(C)\in F(G)$. If both $(C-v_{e})+u$ and $(C-v_{e})+v$ are in $F(G)$ in the case of Proposition \ref{prop:injection1} (iii), then define $\psi(C)$ to be one of them arbitrarily. In all cases, it is easy to see that $\psi(C)/e=C$. This shows the assertion. 
\end{proof}

Now we are ready to prove Lemma \ref{lemma:minor}.

\vspace{0.2cm}
\hspace{-0.43cm}\textbf{Proof of Lemma \ref{lemma:minor}.}
Let $G=(V,E)$ be a 3-connected graph and let $e\in E$ be such that $G/e$ is 3-connected. Let $T_{e}(G)$ be the set of all triangles in $G$ using the edge $e$. It is easy to observe that 
\begin{equation}
|E(G)|-|E(G/e)|=|T_{e}(G)|+1.
\end{equation}
On the other hand, the image \text{im $\psi$} of $\psi$ is disjoint from $T_{e}(G)$. Otherwise, there exists $C\in F(G/e)$ such that $\psi(C)=T$ for some $T\in T_{e}(G)$. But then $C = \psi(C)/e = T/e \simeq K_{2}$, a contradiction. This gives us 
\begin{equation}
|F(G)|-|\text{im $\psi$}|\ge |T_{e}(G)|.
\end{equation}
Then we have 
\begin{align}
\Lambda_{F}(G) - \Lambda_{F}(G/e) &= (|F(G)|-|F(G/e)|) - (|E(G)|-|E(G/e)|) + 1 \\
& =  (|F(G)|-|\text{im $\psi$}|)-(|T_{e}(G)|+1)+1\ge 0,
\end{align}
as desired. This shows the assertion. $\hfill\blacksquare$
\vspace{0.2cm}

We conclude this section by giving a short proof of Tutte's cycle space theorem, using the injection $\psi:F(G/e)\rightarrow F(G)$ that we used in proving Theorem \ref{thm:hadwiger}. Given a simple graph $G$, its \textit{edge space}, denoted $\mathcal{E}(G)$, is a free $\mathbb{Z}_{2}$-module with basis $\{\mathbf{1}(e)\,:\, e\in E(G)\}$, where $\mathbf{1}(e)$ is the indicator function of the edge $e$. For any subgraph $H\subseteq G$, denote $\mathbf{1}(H) := \sum_{e\in E(H)} \mathbf{1}(e)$. The submodule of $\mathcal{E}(G)$ generated by the elements in $\{ \mathbf{1}(C) \,:\, \text{$C\subseteq G$ is a cycle} \}$ is called the \textit{cycle space} of $G$, which we denote by $\mathcal{C}(G)$. Denote by $\mathcal{F}(G)$ the submodule of $\mathcal{C}(G)$ generated by the elements in $\{ \mathbf{1}(F) \,:\, F\in F(G) \}$.

\begin{customthm}{3.5}[Tutte \cite{tutte1963draw}]
	For any 3-connected graph $G$, we have $\mathcal{C}(G)=\mathcal{F}(G)$.  
\end{customthm}

\begin{proof}
	In this proof, all summations are taken modulo 2. We may assume $|V|\ge 5$ since otherwise $G$ is a complete graph and the assertion is clear. Suppose, for a contradiction, that there exists a counterexample $G=(V,E)$ with $|V(G)|$ minimal. Let $e=uv$ be an edge in $G$ such that $G/e$ is 3-connected (the existence of such ``contractible'' edge follows from Theorem \ref{seymour}, but a short proof of which is given in  \cite[Lemma 3.2.1]{diestel2000graph}). Then $\mathcal{C}(G/e)=\mathcal{F}(G/e)$ and $T_{e}(G)\subseteq F(G)$, where $T_{e}(G)$ denotes the set of triangles in $G$ using the edge $e$.

	Note that for any subgraph $H\subseteq G$,  $\mathbf{1}(H)\in \mathcal{C}(G)$ if and only if $H$ is the union of some edge-disjoint cycles in $G$. If we define $T(H)$ to be the set of all vertices $w\in V(H)$ such that $uw,vw\in E(H)$, then note that for any $H\subseteq G$ with $\mathbf{1}(H)\in \mathcal{C}(G)$, $\mathbf{1}(H/e)\in \mathcal{C}(G/e)$ if and only if $T(H)=\emptyset$. 
	
	Choose any subgraph $\alpha\subseteq G$ such that $\mathbf{1}(\alpha)\in \mathcal{C}(G)\setminus \mathcal{F}(G)$. Since $T_{e}(G)\subseteq F(G)$, observe that $\mathbf{1}(\alpha)+\sum_{w\in T(\alpha)} \mathbf{1}(T_{uwv}) \in \mathcal{C}(G)\setminus \mathcal{F}(G)$, where $T_{abc}$ denotes the triangle on vertex set $\{ a,b,c \}$. Also note that if we let $\tilde{\alpha}\subseteq G$ be the subgraph such that $\mathbf{1}(\tilde{\alpha})= \mathbf{1}(\alpha)+\sum_{w\in T(\alpha)}\mathbf{1}(T_{uwv})$ with $|V(\tilde{\alpha})|$ minimal, then $T(\tilde{\alpha})=\emptyset$. Hence we see that the following collection of subgraphs
	\begin{equation}
	\mathcal{B}(G) = \{ \alpha\subseteq G\,:\, \text{$\mathbf{1}(\alpha)\in \mathcal{C}(G)\setminus \mathcal{F}(G)$ and $T(\alpha)=\emptyset$} \}
	\end{equation}
	is nonempty. Since $\mathcal{C}(G/e) = \mathcal{F}(G/e)$, for any $\alpha\in \mathcal{B}(G)$, there exists  $F_{1},\cdots,F_{n}\in F(G/e)$ such that  $\mathbf{1}(\alpha/e)=\sum_{i=1}^{n} \mathbf{1}(F_{i})$. We may choose $\alpha\in \mathcal{B}(G)$ so that $n\ge 1$ is as small as possible.

	Let $\psi:F(G/e)\rightarrow F(G)$ be the injection defined in the proof of Proposition \ref{prop:injection2}. We claim that there exists a (possibly empty) subset $T_{0}\subseteq T_{e}(G)$ such that 
	\begin{equation}\label{eq:tutte_pf_claim}
	\mathbf{1}(\alpha) + \mathbf{1}(\psi(F_{n})) + \sum_{T\in T_{0}} \mathbf{1}(T)  =  \mathbf{1}(\alpha') 
	\end{equation} 
	for some $\alpha'\subseteq G$ such that $T(\alpha')=\emptyset$. To see the assertion follows from the claim, note that (\ref{eq:tutte_pf_claim}) yields $\mathbf{1}(\alpha')\in \mathcal{C}(G)$. Moreover, since $T(\alpha')=\emptyset$  it holds that 
	\begin{eqnarray}
	\mathbf{1}(\alpha'/e) = \mathbf{1}(\alpha/e)+ \mathbf{1}(\psi(F_{n})/e) =\left(\sum_{i=1}^{n}\mathbf{1}(F_{i})\right) + \mathbf{1}(F_{n}) = \sum_{i=1}^{n-1}\mathbf{1}(F_{i}).
	\end{eqnarray}
	Then the minimality of $n$ implies $\mathbf{1}(\alpha')\in \mathcal{F}(G)$. Since $T_{e}(G)\subseteq F(G)$, all indicators in (\ref{eq:tutte_pf_claim}) except $\mathbf{1}(\alpha)$ belong to $\mathcal{F}(G)$. This implies $\mathbf{1}(\alpha)\in \mathcal{F}(G)$, a contradiction.

	To show the claim, first let $\beta\subseteq G$ be such that $\mathbf{1}(\beta) = \mathbf{1}(\alpha)+\mathbf{1}(\psi(F_{n}))$ with $|V(\beta)|$ minimal. If $T(\beta)=\emptyset$, then the claim holds for $T_{0}=\emptyset$. Otherwise, suppose there exists a vertex $w\in T(\beta)$, so $uw,vw\in E(\beta)$. Since $T(\alpha)=\emptyset$, either $uw$ or $vw$ is in $E(\psi(F_{n}))$. So $v_{e}\in V(F_{n})$. Let $x,y\in V(F_{n})$ be the two neighbors of $v_{e}$ in $F_{n}\subseteq G/e$. Then $x,y\in V(G)$ and $z$ must be either $x$ or $y$ since $F_{n}$ is induced in $G/e$. This shows $T(\beta)\subseteq \{ x,y \}$. Without loss of generality, we assume $x\in T(\beta)$. If $T(\beta)=\{ x \}$, then the claim holds with $T_{0} = \{ T_{uxv}  \}$. Otherwise, $T(\beta)=\{x,y\}$ and the claim holds with $T_{0} = \{T_{uxv}, T_{uyv} \}$. This shows the assertion.  
\end{proof}

\vspace{0.2cm}
\section{Concluding remarks}

For a 3-connected planar graph $G$, the parameter $\Lambda_{F}(G)$ agrees with the Euler characteristic, so $\Lambda_{F}(G)=2$. Hence Theorem \ref{thm:hadwiger} implies that $h(G)\le 4$. So if Theorem \ref{thm:chromatic} could be shown for the parameter $\Lambda_{F}$ instead of $\Lambda_{C}$, this implies $\chi(G)\le 4$ for all 3-connected planar $G$. This is equivalent to the famous four-color theorem. Proving such a result would require a ``chain theorem'' for homomorphisms, playing a similar role of the Splitter theorem for edge contractions. Indeed, Kriesell \cite{kriesell1998contractible} showed that 3-connected non-complete graphs always have a contractible nonedge, i.e., a pair of non-adjacent vertices $u,v$ such that identifying them yields a 3-connected graph. However, it seems that cycles behave less nicely under nonedge contractions than under edge contractions, as identifying two nonadjacent vertices could make a non-separating induced path into a non-separating induced cycle. 

In Lemma \ref{lemma:minor} and Theorem \ref{thm:hadwiger}, we have shown that if $K$ is a complete graph minor of a 3-connected graph $G$, then $\Lambda_{F}(K)\le \Lambda_{F}(G)$. Does this hold for any 3-connected minor $H\le_{m} G$? In fact, Lemma \ref{lemma:minor} and a similar argument in the proof of Theorem \ref{thm:hadwiger} shows that $\Lambda_{F}(H)\le \Lambda_{F}(G)$ whenever $H$ is a 3-connected minor of $G$ that can be obtained by a sequence of edge contractions. For general minors involving edge deletions, one would naturally try to show a monotonicity of $\Lambda_{F}$ under edge deletions, which is similar to Lemma \ref{lemma:minor}. Namely, this reads that for every 3-connected non-complete graph $G$, there exists an edge $e$ such that $G-e$ is 3-connected and 
\begin{equation}\label{deletion}
\Lambda_{F}(G-e)\le \Lambda_{F}(G).
\end{equation} 
Note that (\ref{deletion}) is equivalent to $|F(G-e)|\le |F(G)|-1$.

Indeed, Thomassen \cite{thomassen1980planarity} showed that in a 3-connected graph, there are at least two induced non-separating cycles using any fixed edge $e$, just like the facial cycles of a planar graph do. So at least one element of $F(G)$ is deleted in $G-e$. However, it could be that many non-separating cycles in $G$ have ``chord'' $e$ so that deletion of which makes all of them induced. More specifically, let $G$ be a graph with a vertex $w$ such that $w$ is adjacent to all vertices in $H$, where $H$ is a 2-connected graph consisting of $n$ vertex-disjoint induced paths $P_{1},\cdots,P_{n}$ from a vertex $s$ to a vertex $t$. Suppose that $e=st\in E(G)$ and all of $P_{i}$'s have length $\ge 2$. Note that every induced cycle in $H$ and $H-e$ is non-separating in $G$ and $G-e$, respectively. Now $H$ has exactly $n$ induced cycles, whereas $H-e$ has $n(n-1)/2$ of them. Consequently, $|F(G)|=n+|E(H)|$ but $|F(G-e)|=n(n-1)/2+|E(H)|-1$. Therefore $\Lambda_{F}(G- e)>\Lambda_{F}(G)$ when $n\geq  4$.

\vspace{0.2cm}

\section*{Acknowledgment}

The author is grateful for the anonymous referees for value suggestions. The author thanks to Neil Robertson and Paul Seymour for helpful discussions.

	\vspace{0.2cm}

	\bibliographystyle{plain}
	\bibliography{mybib}

\end{document}